\documentclass[leqno, 12pt]{amsart}

\usepackage{amssymb}
\usepackage{amsmath}
\usepackage{enumerate}
\usepackage{amsfonts}
\usepackage{color}

\numberwithin{equation}{section}

\usepackage[colorlinks,citecolor=blue,urlcolor=blue]{hyperref}

\newcommand{\Hla}{\mathcal{H}}

\newcommand{\la}{\lambda}
\newcommand{\supp}{\mathrm{supp}}

\renewcommand{\(}{\left(}
\renewcommand{\)}{\right)}
\newcommand{\tm}{\widetilde{m}}

\newcommand\ssf{\hspace{.25mm}}
\newcommand\ssb{\hspace{-.25mm}}

\newtheorem{remark}[equation]{Remark}
\newtheorem{lem}[equation]{Lemma}
\newtheorem{thm}[equation]{Theorem}
\newtheorem{pro}[equation]{Proposition}

\begin{document}
\title[Multivariate multiplier theorem for the Hankel transform]{Multivariate H\"{o}rmander-type multiplier theorem\\ for the Hankel transform}
\author{Jacek Dziuba\'{n}ski,  Marcin Preisner, and B\l a\.{z}ej Wr\'obel}
\address{Instytut Matematyczny, Uniwersytet Wroc\l awski, pl. Grunwaldzki 2/4, 50-384 Wroc\l aw, Poland}
\email{\ \newline
jdziuban@math.uni.wroc.pl,\newline
preisner@math.uni.wroc.pl (corresponding author), \newline
blazej.wrobel@math.uni.wroc.pl}
\thanks{The research was partially supported by Polish funds for sciences, grants: N N201 397137 and N N201 412639, MNiSW, and research project 2011/01/N/ST1/01785, NCN}

\begin{abstract}
Let $\Hla (f)(x)=\int_{(0,\infty)^d} f(\la) E_{x}(\la)\,d\nu(\la),$ be the multivariate Hankel transform,
where $E_{x}(\la)=\prod_{k=1}^d (x_k \la_k)^{-\alpha_k+1/2}J_{\alpha_k-1/2}(x_k \la_k)$, with $d\nu(\lambda)=\lambda^\alpha \, d\lambda$, $\alpha=(\alpha_1,...,\alpha_d)$. We give sufficient conditions on a bounded function $m(\lambda)$ which guarantee that the operator $\mathcal H(m\mathcal H f)$ is bounded on $L^p(d\nu)$ and of  weak-type (1,1), or bounded on the Hardy space $H^1((0,\infty)^d,d\nu)$ in the sense of Coifman-Weiss.
\end{abstract}

\maketitle

\section{Introduction and preliminaries}
For a multiindex $\alpha =(\alpha_1,...,\alpha_d)$, $\alpha_k>-1\slash 2$, we consider the measure space $X=((0,\infty)^d, d\nu(x))$, where  $d\nu(x)=d\nu_{1}(x_1)\cdots d\nu_{d}(x_d),$ $d\nu_{k}(x_k)=x_k^{2\alpha_k}\,dx_k,$ $k=1,\ldots,d.$
The space $X$ equipped with the Euclidean distance is a space of homogeneous type in the sense of Coifman-Weiss. We denote by $H^1(X)$ the atomic Hardy space associated with $X$ in the sense of \cite{CW}. More precisely, we say that a measurable function $a$ is an $H^1(X)$-atom, if there exists a ball $B$, such that $\supp\, a \subset B,$ $\|a\|_{L^{\infty}(X)}\leq 1/ \nu(B),$ and $\int_{(0,\infty)^d}a(x)d\nu(x) =0.$ The space $H^1(X)$ is defined as the set of all $f\in L^1(X),$ which can be written as $f= \sum_{j=1}^{\infty} c_j a_j,$ where $a_j$ are atoms and $\sum_{j=1}^{\infty} |c_j|<\infty,$ $c_j\in\mathbb{C}.$ We equip $H^1(X)$ with a norm
\begin{equation}
\label{harnorm}
\|f\|_{H^1(X)}=\inf \sum_{j=1}^{\infty} |c_j|,
\end{equation}
where the infimum is taken over all absolutely summable sequences $\{c_j\}_{j\in\mathbb{N}},$ for which $f= \sum_{j=1}^{\infty} c_j a_j,$ with $a_j$ being $H^1(X)$-atoms.

For an appropriate function $f$ the (modified) Hankel transform is defined by

\begin{equation*}
\Hla (f)(x)=\int_{(0,\infty)^d} f(\la) E_{x}(\la)\,d\nu(\la),
\end{equation*}
where $$E_{x}(\la)=\prod_{k=1}^d (x_k \la_k)^{-\alpha_k+1/2}J_{\alpha_k-1/2}(x_k \la_k)=\prod_{k=1}^dE_{x_k}(\la_k).$$ Here $J_{\nu}$ is the Bessel function of the first kind of order $\nu,$ see \cite[Chapter 5]{3}. The system $\{E_{x}\}_{x\in(0,\infty)^d}$ consists of the eigenvectors of the Bessel operator
\begin{equation*}
L=-\Delta-\sum_{k=1}^d \frac{2\alpha_k}{\la_k}\frac{\partial}{\partial \la_k};
\end{equation*}
that is, $L(E_{x})=|x|^2E_{x}.$ Also, the functions $E_{x_k},$ $k=1,\ldots,d$, are eigenfunctions of the one-dimensional Bessel operators
$$L_{k}=-\frac{\partial^2}{\partial {\la_k}^2}-\frac{2\alpha_k}{\la_k}\frac{\partial}{\partial \la_k},$$
namely, $L_{k}(E_{x_k})=x_k^2 E_{x_k}.$

It is known that $\Hla$ is an isometry on $L^2(X)$ that satisfies $\Hla^{-1}=\Hla$ (see, e.g., \cite[Chapter 8]{T}).  Moreover, for  $f\in L^2(X),$ we have \begin{equation}\label{diago}L_k(f)=\Hla(\la_k^2\Hla f).\end{equation}

For $y\in X$ let $\tau^{y}$ be the $d$-dimensional generalized Hankel translation given by
$$ \Hla (\tau^yf)(x)=E_y(x)\Hla f(x).$$
Clearly,
$\tau^y f (x)=\tau^{y_1}\cdots\tau^{y_d}f (x)$, where for each  $k=1,\ldots,d$, the operator  $\tau^{y_k}$  is the one-dimensional Hankel translation acting on a function $f$ as a function of the  $x_k$  variable with the other variables fixed.  It is also known that $\tau^y$ is a contraction on all $L^p(X)$ spaces, $1\leq p\leq\infty,$ and that
$$ \tau^yf(x)=\tau^xf(y).$$
For two reasonable functions $f$ and $g$ define their Hankel convolution as
$$f \natural g(x)= \int_{X}\, \tau^{x}f(y) g(y)\,d\nu(y).$$
It is not hard to check that $f\natural g=g\natural f$ and
\begin{equation}\label{convo}\Hla(f\natural g)(x)=\Hla f(x)\,\Hla g(x).\end{equation}
As a consequence of the contractivity of $\tau^{y}$ we also have
\begin{equation}
\label{young}
\|f\natural g\|_{L^1(X)}\leq \|f\|_{L^1(X)} \|g\|_{L^1(X)}, \qquad f\in L^1(X),\quad g\in L^1(X).
\end{equation}

For details concerning translation, convolution, and transform in the Hankel setting we refer the reader to, e.g., \cite{Ha}, \cite{T}, and \cite{W}.

For a function $f\in L^1(X)$ and $t>0$ let $f_t$ denote the $L^1(X)$-dilation of $f$ given by $$(f_t)(x)={t}^{Q}f(tx),$$ where $Q=\sum_{k=1}^d (2\alpha_k+1)$. Then we have:
\begin{align}
&\Hla (f_t)(x)=\Hla f(t^{-1} x)\label{dil},
\\&\label{diltra}\tau^{y}(f_t)(x)=(\tau^{ty}f)_t(x).
\end{align}
Notice that $Q$ represents the dimension of $X$ at infinity, that is, $\nu(B(x,r))\sim r^Q$ for large $r$.

 Let $m:X\rightarrow \mathbb{C}$ be a bounded measurable function. Define the multiplier operator $\mathcal T_m$ by
\begin{equation}
\label{spHan}
\mathcal T_m(f)=\Hla(m\Hla f).
\end{equation}
Clearly, $\mathcal T_m$ is bounded on $L^2(X).$ Also note that if $m(\la_1,\ldots,\la_d)=n(\la_1^2,\ldots,\la_d^2),$ for some bounded, measurable function $n$ on $\mathbb{R}^d,$ then from \eqref{diago} it can be deduced that the Hankel multiplier operator defined by \eqref{spHan} coincides with the joint spectral multiplier operator $n(L_1,\ldots,L_d).$ The smoothness requirements on $m$ that guarantee the boundedness of $\mathcal T_m$ on, e.g., $L^p(X)$ will be stated in terms of appropriate Sobolev space norms.

For $z\in \mathbb C$, $\text{Re}\, z>0$, let
$$ G_z(x)=\Gamma \big(z\slash 2\big)^{-1} \int_0^\infty (4\pi t)^{-d\slash 2 }e^{-|x|^2\slash 4t} e^{-t}t^{z\slash 2}\frac{dt}{t}$$
be the kernels of the Bessel potentials. Then
\begin{equation}\label{potential}\| G_z\|_{L^1(\mathbb R^d)}\leq \Gamma (\text{Re} \, z\slash 2)|\Gamma (z\slash 2)|^{-1} \ \ \text{ and}\ \  \mathcal F G_z(\xi)=(1+|\xi|^2)^{-z\slash 2},
\end{equation}
where $\mathcal F G_z(\xi)=\int_{\mathbb{R}^d}G_z(x)e^{-i<x,\xi>}\,dx$ is the Fourier transform.

By definition, a function $f\in W_2^s(\mathbb R^d)$, $s>0$,  if and only if there exists a function $h\in L^2(\mathbb R^d)$ such that $f=h\star G_s$, and $\| f\|_{W^s_2(\mathbb R^d) } = \| h\|_{L^2(\mathbb R^d)}$.

 Similarly, a function $f$ belongs to the potential space $\mathcal L_s^\infty(\mathbb R^d)$, $s>0$,  if there is a function $h\in L^\infty(\mathbb R^d)$ such that $f=h\star G_s$ (see \cite[Chapter V]{Stein2}). Then $\| f\|_{\mathcal L^\infty_s (\mathbb R^d)}=\| h\|_{L^\infty (\mathbb R^d)}$.

Denote $A_{r,R} =\{x\in\mathbb{R}^d\,:\, r\leq|x|\leq R\}.$ The main results of the paper are the following theorems.
\begin{thm}
\label{Thmoo} Assume that $\alpha_k\geq 1\slash 2$ for $k=1,...,d$.
 Let $m(\la)=n(\la_1^2,\ldots,\la_d^2),$ where $n$ is a bounded function on $\mathbb{R}^d$ such that, for certain real number $\beta > Q/2$ and for some (equivalently, for every) non-zero radial function $\eta\in C_c^{\infty}(A_{1/2,2}),$ we have
\begin{equation}
\label{conmulto}
\sup_{j\in\mathbb{Z}}\|\eta(\cdot)n(2^j\cdot)\|_{W^{\beta}_2(\mathbb R^d)}\leq C_{\eta}.
\end{equation}
Then the multiplier operator $\mathcal T_m$ is a Calder\'{o}n-Zygmund operator associated with the kernel
$$K(x,y)=\sum_{j\in\mathbb{Z}}\tau^y \Hla(\psi(2^{-j}(\la_1^2,\cdots,\la_d^2)) m(\la)) (x),$$
where $\psi$ is a $C_c^{\infty}(A_{1/2,2})$ function such that
\begin{equation}\label{suma}
\sum_{j\in\mathbb{Z}} \psi(2^{-j}\la)=1, \qquad \lambda\in \mathbb R^d\backslash \{0\}.
\end{equation}
As a consequence $\mathcal T_m$ extends to the bounded operator from $L^1(X)$ to $L^{1,\infty}(X)$ and from $L^p(X)$ to itself for  $1<p<\infty$.
\end{thm}

\begin{thm}
\label{Thmhar} Assume that $\alpha_k\geq 1\slash 2$ for $k=1,...,d$.
 Let $m(\la)=n(\la_1^2,\ldots,\la_d^2),$ where $n$ is a bounded function on $\mathbb{R}^d$ such that, for certain real number $\beta > Q/2$ and for some (equivalently, for every) non-zero radial function $\eta\in C_c^{\infty}(A_{1/2,2}),$ \eqref{conmulto} holds.
Then the multiplier operator $\mathcal T_m$ extends to a bounded operator on the Hardy space $H^1(X).$
\end{thm}

\begin{remark}\label{remark1}
 If we relax the conditions on $\alpha_k$ assuming only that $\alpha_k> -1\slash 2$, then the conclusions of
 Theorems \ref{Thmoo} and \ref{Thmhar} hold provided there is $\beta >Q\slash 2 $ such that
 \begin{equation}\label{conmulto22}
\sup_{j\in\mathbb{Z}}\|\eta(\cdot)n(2^j\cdot)\|_{\mathcal L_{\beta}^\infty(\mathbb R^d)}\leq C_{\eta}.
\end{equation}
\end{remark}

The weak type $(1,1)$ estimate under assumption \eqref{conmulto22} could be proved by applying a general multiplier theorem of Sikora \cite{S}. However, in the case of the  Hankel transform Remark \ref{remark1} has a simpler proof based on Lemma \ref{Lporr} and Remark \ref{remark2}.

 Hankel multipliers, mostly of one variable,  attracted attention of many authors, see, e.g.,   \cite{BBB}, \cite{2}, \cite{1}, \cite{GS}, \cite{GT}, and references therein. In \cite{BBB} the authors considered multidimensional Hankel multipliers $m$ of Laplace transform type, that is,
 $$ m(y)=|y|^2\int_0^{\infty} e^{-t|y|^2}\phi(t)\, dt,$$
 where $\phi \in L^{\infty} (0,\infty)$ (see \cite{Stein}). Setting
 $$n(\lambda)=\Xi(\lambda) (\lambda_1+\ldots+\lambda_d) \int_0^\infty e^{-t(\lambda_1+...+\lambda_d)}\phi (t)\, dt,$$ where $\Xi \in C^\infty (\mathbb R^d\setminus \{0\})$, $\Xi (t\lambda)=\Xi (\lambda)$ for $t>0$, $\Xi (\lambda )=1$ for $\lambda \in (0,\infty)^d$, $\Xi (\lambda)=0$ for $\lambda_1+...+\lambda_d<|\lambda |\slash d,$ we easily see that (\ref{conmulto}) and (\ref{conmulto22}) hold with every $\beta>0$.

 For other results and references concerning spectral multiplier theorems on $L^p$ spaces the reader is referred to
 \cite{A},  \cite{C}, \cite{He}, \cite{MM}, \cite{Ma}, \cite{MS}, and \cite{S}.

\vskip ,5em

{\bf Acknowledgments.} The authors would like to thank Alessio Martini for discussions on spectral multipliers.

\section{Auxiliary estimates}
In this section we prove some basic estimates needed in the sequel. Denote $w^s(x)=(1+|x|)^s$.
\begin{lem}
\label{Lporr}
 For every $s,\varepsilon>0$ there exists a constant $C_{s,\varepsilon}$ such that   if $m(\la)=n(\la_1^2,...,\la_d^2)$, $\supp \,n\subseteq A_{1/4,4}$, then
\begin{equation}\label{eq21}
\|\Hla (m) w^s\|_{L^2(X)} \leq C_{s,\varepsilon} \|n\|_{W^{s+d\slash 2+\varepsilon}_2(\mathbb R^d)}.
\end{equation}
\end{lem}

{\bf Proof.} Since $m(\la)= g(\la_1^2,...,\la_d^2) e^{-|\la|^2},$ with $g(\la)= n(\la)e^{\la_1+...+\la_d},$ using the Fourier inversion formula for $g$, we get
\begin{equation*}\begin{split} (2\pi)^d\, m(\la)&=e^{-|\la|^2} \int_{\mathbb R^d}\mathcal{F}(g)(y)e^{i y_1\la_1^2+...+iy_d\la_d^2}\,dy\\
&=\int_{\mathbb R^d} \mathcal{F}(g)(y)e^{(-1+iy_1)\la_1^2+...+(-1+iy_d)\la_d^2}\,dy.
     \end{split}\end{equation*}
 Applying the Hankel transform and changing the order of integration, we obtain
\begin{equation}\label{inte}\Hla (m)(x)=(2\pi)^{-d}\int_{\mathbb R^d} \mathcal{F}(g)(y)\Hla (e_{\textbf{1}-iy})(x)\,dy,\end{equation}
where for $z=(z_1,\ldots,z_d)\in\mathbb{C}^d,$ $e_z(\la)=\prod_{k=1}^d e_{z_k}(\la_k)$ with $e_{z_k}(\la_k)=e^{-z_k\la_k^2},$ while $\textbf{1}=(1,\ldots,1).$ Clearly, $$\Hla(e_{\textbf{1}-iy})(x)=\prod_{k=1}^{d}\Hla_k(e_{1-iy_k})(x_k),$$ with $\Hla_k$ denoting the one-dimensional Hankel transform acting on the $k$-th variable.
 It is well known that for $t>0,$ $\Hla_k(e_{t})(x_k)=Ct^{-(2\alpha_k+1)/2}\exp{(-x_k^2\slash 4t)},$ see \cite[p. 132]{3}. Moreover, for  fixed $x_k$, the functions $$z_k\mapsto \Hla_k(e_{z_k})(x_k)\qquad \textrm{and} \qquad z_k\mapsto Cz_k^{-(2\alpha_k+1)/2}\exp{\left(-\frac{x_k^2}{4{z_k}}\right)}$$ are holomorphic on $\{z_k\in\mathbb{C}\,:\, \textrm{Re} \, z_k>0\}$ (provided we choose   an appropriate holomorphic branch of the  power function $z_k^{-(2\alpha_k+1)/2}$). Hence, by the uniqueness of the holomorphic extension, we obtain $$\Hla_k(e_{1-iy_k})(x_k)=C (1-iy_k)^{-(2\alpha_k+1)/2}\exp{\bigg(-\frac{x_k^2}{4(1-iy_k)}\bigg)}.$$
Since $\textrm{Re}\, x_k^2\slash 4(1-iy_k)=x_k^2\slash 4(1+y_k^2),$ the change of variable $x_k=(1+y_k^2)^{1/2}u_k$  leads to \begin{equation}\label{onprod}\int_{(0,\infty)}|x_k^s\Hla(e_{1-iy_k})(x_k)|^2\,d\nu_k(x_k)\lesssim (1+y_k^2)^{s},\qquad s\geq0.\end{equation}
Now, observing that $(1+|x|)^{2s}\approx 1+x_1^{2s}+\ldots +x_d^{2s}$ and using \eqref{onprod} we arrive at
$$\|(1+|\cdot|)^s\Hla(e_{\textbf{1}-iy})(\cdot )\|_{L^2(X)}\lesssim \sum_{k=1}^d (1+y_k^2)^{s/2}\approx (1+|y|)^{s}.$$
The latter bound together with \eqref{inte}, Minkowski's integral inequality, and the Schwarz  inequality give
\begin{align*}\|\Hla (m) w^s\|_{L^2(X)}&\lesssim \int_{\mathbb{R}^d}|\mathcal{F}(g)(y)|(1+|y|)^{s}\,dy \\&\lesssim \bigg( \int_{\mathbb{R}^d}|\mathcal{F}(g)(y)|^2(1+|y|)^{2s+d+2\varepsilon}\,dy\bigg)^{1/2}\bigg(\int_{\mathbb{R}^d}(1+|y|)^{-d-2\varepsilon}\,dy\bigg)^{1/2}\\
&\lesssim \|g\|_{W^{s+d/2+\varepsilon}_2(\mathbb R^d)}\end{align*}
for any fixed $\varepsilon>0.$
Since $g(\la)=n(\la)e^{\la_1+...+\la_d}=n(\la) (e^{\lambda_1+...+\la_d}\eta_0(\la)),$ for some $\eta_0\in C_c^{\infty}(A_{1/8,8}),$ we see that $\|g\|_{W^{s+d/2+\varepsilon}_2(\mathbb R^d)}\leq C \|n\|_{W^{s+d/2+\varepsilon}_2(\mathbb R^d)},$ which implies \eqref{eq21}. $\square$\\
\indent Remark that a slight modification of the reasoning above shows that if $m(\la)=n(\la_1^2,\ldots,\la_d^2),$ $n\in C_c^{\infty}(A_{1/2,2}),$ then
\begin{equation}
\label{cksob}
|\Hla(m)(x)|\leq C_N \|n\|_{C^{N+d}(A_{1/2,2})} w^{-N}(x),
\end{equation}
where $C^N$ denotes the supremum norm on the space of $N$-times continuously differentiable functions.

 Using ideas of Mauceri-Meda \cite{MM} combined with  the fact that the Hankel transform is  an $L^2$-isometry we can improve Lemma \ref{Lporr} in the following way.
\begin{lem}
\label{Linterr}
Assume that $\alpha_k\geq 1\slash 2$ for $k=1,...,d$. Then  for every  $s,\varepsilon>0,$ there is a constant $C_{s,\varepsilon}$ such that
if  $m(\la)=n(\la_1^2, ...,\la_d^2)$,  $\supp \, n\subseteq A_{1/2,2}$, then
\begin{equation*}
\|\Hla (m) w^s\|_{L^2(X)} \leq C_{s,\varepsilon} \|n\|_{W_2^{s+\varepsilon}(\mathbb R^d)}.
\end{equation*}
\end{lem}

\textbf{Proof.}
Let $h\in L^2(\mathbb R^d)$ be such that $n=h\star G_{s+\varepsilon}$. Set $s'=(s+\varepsilon)(d+6)\slash 2\varepsilon$,
\newline $\theta=2\varepsilon \slash  (6+d)$.
Define $n_z$ by
 \begin{equation*}
\mathcal{F}(n_{z})(\xi)=\mathcal F h(\xi)(1+|\xi|^2)^{-s'z\slash 2}, \ \ 0\leq  \text{Re}\, z\leq 1.
\end{equation*}
 Clearly, $n_z=h\star G_{s'z}$,  $\text{Re}\, z>0$, and $n=n_\theta$.
 Let $\eta_0$ be a $C_{c}^{\infty}$ function supported in $A_{1/4,4},$ equal to $1$ on $A_{1/2,2},$ and let $N_z(\la)=n_z(\la)\eta_0(\la)$. Then $\supp \, N_z\subseteq A_{1/4,4}$ and $\mathcal{F}(N_z)=\mathcal{F}(n_z)\star \mathcal{F}(\eta_0).$
  Define
  $$m_z(\la)=n_z(\la_1^2,...,\la_d^2)\quad \textrm{ and } \quad M_z(\la)=N_z(\la_1^2,...,\la_d^2).$$
Since $\alpha_k\geq 1\slash 2$  for every $k=1,...,d$, we have that $M_z \in L^2(X)$ and $\|M_z\|_{L^2(X)} \lesssim \|N_z\|_{L^2(\mathbb R^d)}$.
Let $g$ be an arbitrary $C_c^{\infty}(X)$ function with $\|g\|_{L^2(X)}=1.$ Set
\begin{equation}\label{functionF}F(z)=\int_{X} \Hla(M_z)(x)(1+|x|)^{(s'-3-d\slash 2)z}g(x)\,d\nu(x).
\end{equation}
Then $F$ is holomorphic in the strip $S=\{z: 0<\textrm{Re}\, z<1\}$
and also continuous and bounded on its closure $\bar S$.
Using Parseval's equality and the facts that $\supp \,N_z\subseteq A_{1/4,4}$ and $\mathcal{F}(\eta_0)\in\mathcal{S}(\mathbb{R}^d),$ for $\textrm{Re} \,z=0,$ we get
\begin{align*}
|F(z)|&\leq \|\Hla(M_z)\|_{L^2(X)}=\|M_z\|_{L^2(X)}\leq C \|N_z\|_{L^2(\mathbb{R}^d)}\approx \|\mathcal{F}N_z\|_{L^2(\mathbb{R}^d)}\\
&\leq C_{\eta_0,s',\theta} \|n\|_{W_2^{s+\varepsilon}(\mathbb R^d)}.
\end{align*}
If $\textrm{Re}\, z=1$, then  applying in addition Lemma \ref{Lporr},   we obtain
\begin{align*}
|F(z)|& \leq \|\Hla(M_z)w^{s'-3-d\slash 2}\|_{L^2(X)}  \leq C\|N_z\|_{W_2^{s'}(\mathbb R^d)}\\
&\leq C_{\eta_0}\|n_z\|_{W_2^{s'}(\mathbb R^d)}=C\| h\|_{L^2(\mathbb R^d)} = C \|n\|_{W_2^{s+\varepsilon}(\mathbb R^d)}.
\end{align*}
From the Phragm\'{e}n-Lindel\"{o}f principle we get  $|F(\theta)|\leq C\|n\|_{W^{s+\varepsilon}_2(\mathbb R^d)}.$ Taking the supremum over all such $g$ we arrive at
$$\|\Hla(M_{\theta})w^{(s'-3-d/2)\theta}\|_{L^2(X)}\leq C \|n\|_{W_2^{s+\varepsilon}(\mathbb R^d)}.$$
Recall that $n=n_\theta=N_\theta$, so that also $m=m_\theta = M_{\theta},$ hence we get the desired conclusion. $\square$

\begin{remark}\label{remark2}
If we relax the conditions on $\alpha_k$ in Lemma \ref{Linterr} by assuming that $\alpha_k > -\frac{1}{2}$, then
\begin{equation*}
\|\Hla (m) w^s\|_{L^2(X)} \leq C_{s,\varepsilon} \|n\|_{\mathcal L^\infty_{s+\varepsilon}(\mathbb R^d)}.
\end{equation*}
\end{remark}

{\bf Proof.} We argue similarly to the proof of Lemma \ref{Linterr}. Indeed, write $n=h\star G_{s+\varepsilon}$, where $h\in L^\infty (\mathbb R^d)$.
Since $\supp \, n\subset A_{1\slash 2, 2}$, one can prove that $h\in L^2(\mathbb R^d)$ and $\| h\|_{L^2(\mathbb R^d)}\leq C_{s,\varepsilon} \| n\|_{\mathcal L^\infty_{s+\varepsilon}}$.

Set $s'=(2s+\varepsilon)(6+d)\slash 2\varepsilon$, $\theta=\varepsilon\slash (6+d)$ and
define
$$
N_z(\lambda) =\eta_0(\lambda)\, h\star G_{s'z+\varepsilon \slash 2}(\lambda), \ \ \lambda\in\mathbb R^d, \ \ 0\leq \text{Re}\, z\leq 1.$$
Then for every $z\in \bar S$ the function $N_z(\lambda) $ is continuous and supported in $A_{1\slash 4,4}$. Let $M_z(\lambda)=N_z(\lambda_1^2,...,\lambda_d^2)$.  Clearly, $M_\theta=m$. Moreover, by (\ref{potential}),
$$\|M_z \|_{L^2(X)}\leq C\| M_z\|_{L^\infty(X)}\leq C\| N_z\|_{L^\infty(\mathbb R^d)}\leq C_{s,\varepsilon} \| h\|_{L^\infty(\mathbb R^d)}=C_{s,\varepsilon}\|n\|_{\mathcal L^\infty_{s+\varepsilon}(\mathbb R^d)}.$$
We now use the new functions $M_z$ to define a bounded  holomorphic  function $F(z)$ by the formula (\ref{functionF}).
Obviously $|F(z)|\leq C_{s,\varepsilon} \| n\|_{\mathcal L^\infty_{s+\varepsilon}}$ for $\text{Re}\, z=0$. To estimate $F(z)$ for $\text{Re}\, z=1$ we utilize  Lemma \ref{Lporr} and obtain
\begin{equation*}\begin{split} |F(z)|&\leq  \|\Hla(M_z)w^{s'-3-d\slash 2}\|_{L^2(X)}\leq
C\|N_z\|_{W_2^{s'}(\mathbb R^d)} \\
&\leq C_{\eta_0,s,\varepsilon} \| h\|_{L^2(\mathbb R^d)}
\leq  C_{s,\varepsilon} \|n\|_{\mathcal L^\infty_{s+\varepsilon}(\mathbb R^d)}.
\end{split}\end{equation*}
An application of the  Phragm\'{e}n-Lindel\"{o}f principle for $z=\theta$ finishes the proof. \ $\square$
\vskip 1em

 We will also need the following off-diagonal estimate (see \cite[Lemma 2.7]{2}).
\begin{lem}
\label{Lfarr}
Let $\delta>0.$ Then there is $C>0$ such that  for every $y\in X$ and $r,t>0,$ we have
$$
\int_{|x-y|>r}|\tau^{y}(f_t)(x)|\,d\nu(x)\leq C(rt)^{-\delta}\|f\|_{L^1(X, w^{\delta}(x)d\nu(x))}.
$$
\end{lem}

\textbf{Proof.} Let $B$ be the left-hand side of the inequality from the lemma. If $|x-y|>r$ then there is $k\in \{1,...,d\}$ such that  $|x_k-y_k|>r/\sqrt{d}$. Hence,
$$
B\leq \sum_{k=1}^d \int_{|x_k-y_k|>r/\sqrt d}|\tau^{y}(f_t)(x)|\,d\nu(x)
=\sum_{k=1}^dB_k.
$$
It is known that the generalized translations can be also expressed as
\begin{equation}
\label{translaform}
\tau^y f(x)= \int_{|x_1-y_1|}^{x_1+y_1}...\int_{|x_d-y_d|}^{x_d+y_d}\,f(z_1,...,z_d)
\,dW_{x_1,y_1}(z_1)...dW_{x_d,y_d}(z_d),
\end{equation}
with $W_{x_k,y_k}$ being a probability measure supported in $[|x_k-y_k|,x_k+y_k]$ (see \cite{Ha}).  Thus,
\begin{equation*}
B_k=\int_{|x_k-y_k|>r/\sqrt{d}}\Big|\,\int_{|x_1-y_1|}^{x_1+y_1}...
\int_{|x_d-y_d|}^{x_d+y_d}\,(f_t)(z_1,...,z_d)\,dW_{x_1,y_1}(z_1)...dW_{x_d,y_d}(z_d)\Big|\,d\nu(x).
\end{equation*}
Introducing the factor $(z_k t)^{\delta}(z_k t)^{-\delta}$ to  the inner integral in the above formula and denoting $g(x)=|f(x)|x_k^{\delta}$, we see that
\begin{align*}B_k&\leq C (rt)^{-\delta} \int_{X}\int_{|x_1-y_1|}^{x_1+y_1}...
\int_{|x_d-y_d|}^{x_d+y_d}g_t(z)\,dW_{x_1,y_1}(z_1)...dW_{x_d,y_d}(z_d)\,d\nu(x)\\
&\leq C (rt)^{-\delta} \|\tau^y g_t\|_{L^1(X)}\leq C (rt)^{-\delta} \|f\|_{L^1(X, w^{\delta}d\nu)},\end{align*}
where in the last inequality we have used the fact that $\tau^y$ is a contraction on $L^1(X).$ $\square$

Let $T_t(x,y)=\tau^{y}\Hla(e^{-t|\la|^2})(x)$ be the integral kernels of the heat semigroup corresponding to $L$. Clearly,
$$T_t(x,y)=T_t^{(1)}(x_1,y_1)...T_t^{(d)}(x_d,y_d), $$
where $T_t^{(k)}(x_k,y_k)$ is the one-dimensional heat kernel associated with the operator $L_k$.
\begin{lem}
\label{Lneaa}
There is a constant $C >0$ such that
$$
\int_{X}|T_1(x,y)-T_1(x,y')|\,d\nu(x)\leq C|y-y'|,\qquad y,y'\in X.
$$
\end{lem}

\textbf{Proof.} The proof is a direct consequence of the   one-dimensional result, see \cite[Theorem 2.1]{1}, together with the equality  $$\int_0^\infty\,T_1^{(k)}(x_k,y_k)d\nu_k(x_k)=1,\qquad k=1,2,...,d.\ \ \square$$
In the proof of Theorem \ref{Thmhar} the following version of \cite[Lemma 2.5]{2} will be used.
\begin{lem}
\label{Lemhar}
Assume that $f,g \in L^1((0,\infty)^d,w^{\delta}\,d\nu),$ with certain $\delta>0.$ Then:
$$
\|f\natural g\|_{L^1((0,\infty)^d,w^{\delta}\,d\nu)}\leq \|f\|_{L^1((0,\infty)^d,w^{\delta}\,d\nu)} \|g\|_{L^1((0,\infty)^d,w^{\delta}\,d\nu)}.
$$
\end{lem}

\textbf{Proof.} After recalling the representation \eqref{translaform} the proof is analogous to the proof of \cite[Lemma 2.5]{2}. $\square$

\section{Proof of Theorem \ref{Thmoo}}
Assume that \eqref{conmulto} holds for some $\beta> Q /2.$ Fix  $\psi\in C_c^{\infty}(A_{1/2,2})$ satisfying \eqref{suma}. Let $$K(x,y)=\sum_{j\in\mathbb{Z}}K_j(x,y)=\sum_{j\in\mathbb{Z}}\tau^{y}\Hla(m_{j})(x),$$ where $m_j(\la)=\psi(2^{-j}(\la_1^2,...,\la_d^2))m(\la)=(\psi(2^{-j}\cdot)n(\cdot))(\la_1^2,...,\la_d^2).$ To prove that $\mathcal T_m$ is indeed a Calder\'{o}n-Zygmund operator associated with the kernel $K(x,y)$ we need to verify that it satisfies the H\"{o}rmander integral condition, i.e.,
\begin{equation}
\label{horm}
\int_{|x-y|>2|y-y'|}|K(x,y)-K(x,y')|\,d\nu(x)\leq C
\end{equation}
for $y,y' \in X,$ and the association condition
\begin{equation}
\label{asoc}
\mathcal T_m f(x)=\int_{X} K(x,y)f(y)\,d\nu(y)
\end{equation}
for compactly supported $f\in L^{\infty}(X)$ such that $x\notin \supp \, f$.
We start by proving \eqref{horm}.
It suffices to show that
\begin{equation*}
D_j(y,y')=\int_{|x-y|>2|y-y'|}|K_j(x,y)-K_j(x,y')|\,d\nu(x)\leq C_j, \
\text{with} \ \  \sum_{j\in\mathbb{Z}}C_j<\infty.
\end{equation*}
\indent Let $r=2|y-y'|$ and assume first $j>- 2\log_2 r.$
 Let
$$\tm_j(\la)=m_j(2^{j/2} \la)=(\psi(\cdot)n(2^{j}\cdot))(\la_1^2,...,\la_d^2).$$
Note that $\supp \,(\psi(\cdot)n(2^{j}\cdot))\subseteq A_{1/2,2}.$ From \eqref{dil} we see that $$\Hla(m_{j})(x)=2^{jQ/2}\Hla (\tm_j)(2^{j/2} x)=(\Hla(\tm_j))_{2^{j/2}}(x).$$ From the Schwarz inequality, Lemma \ref{Linterr}, and the assumption \eqref{conmulto} we get
\begin{equation}\begin{split}\label{weight}
\int_{X} |\Hla(\tm_j)|w^{\delta}\,d\nu &\leq \bigg(\int_{X}|\Hla(\tm_j)|^2 w^{Q+4\delta}\,d\nu\bigg)^{1/2}\bigg(\int_{X}w^{-Q-2\delta}\,d\nu\bigg)^{1/2}\\
&\leq C_{\delta}\|\psi(\cdot) n(2^{j}\cdot)\|_{W^{\beta}_2(\mathbb R^d)}\leq C_{\delta},
\end{split}\end{equation}
for sufficiently small $\delta>0.$ Consequently, from Lemma \ref{Lfarr} it follows that
\begin{align*}
D_j(y,y')&\lesssim \int_{|x-y|>r} |\tau^y (\Hla(\tm_j))_{2^{j/2}}(x)|\,d\nu(x)+\int_{|x-y'|>r/2} |\tau^{y'} (\Hla(\tm_j))_{2^{j/2}}(x)|\,d\nu(x)\\
&\lesssim (2^{j/2}r)^{-\delta}\int_{X} |\Hla(\tm_j)|w^{\delta}\,d\nu\leq C_{\delta}(2^{j/2}r)^{-\delta},
\end{align*}
so that
$\sum_{j> -2\log_2 r} D_j(y,y')\leq C.$

Assume now $j\leq- 2\log_2 r.$ Decompose $\tm_j(\la)=\tilde{\theta}_j(\la)  e^{-|\la|^2},$ so that we have $\tilde{\theta}_j(\la)=(\psi(\cdot)\exp(\cdot_1+...+\cdot_d)n(2^{j}\cdot))(\la_1^2,...,\la_d^2).$ Clearly, $\psi(\la)e^{\la_1+...+\la_d}$ is a $C_c^{\infty}$ function supported in $A_{1/2,2}.$ Denote $\tilde{\Theta}_j(x)=\Hla(\tilde{\theta}_j)(x).$ Since $\Hla(m_j)=(\Hla(\tm_j))_{2^{j/2}}$ and $\Hla(\tm_j)=\tilde{\Theta}_j \natural \Hla(e^{-|\la|^2})$ (which is a consequence of \eqref{convo}), by using \eqref{diltra}, we get
\begin{align*}
K_j(x,y)-K_j(x,y')&=(\tau^{2^{j/2} y} \Hla(\tm_j))_{2^{j/2}}(x)-(\tau^{2^{j/2} y'} \Hla(\tm_j))_{2^{j/2}}(x)\\
&=\(\tilde{\Theta}_j \natural \(T_1(\cdot,2^{j/2}y)-T_1(\cdot,2^{j/2}y')\)\)_{2^{j/2}}(x).
\end{align*}
Proving \eqref{weight} with $\tm_j$ replaced by
$\tilde{\theta}_j$ and $\delta=0$ poses no difficulty. Hence, from Lemma \ref{Lneaa} and \eqref{young} we obtain
$$D_j(y,y')\leq \|\tilde{\Theta}_j\|_{L^1 (X)}\|T_1(\cdot,2^{j/2}y)-T_1(\cdot,2^{j/2}y')\|_{L^1(X)}\leq C 2^{j/2}|y-y'|.$$
Consequently, $\sum_{j\leq- 2\log_2 r}D_j(y,y')\leq C$ and the proof of \eqref{horm} is finished.\\

 Now we turn to the proof of \eqref{asoc}. From the assumptions, for some $R>r>0,$ $$\int_{X} K_j(x,y)f(y)\,d\nu(y)=\int_{R>|x-y|>r}K_j(x,y)f(y)\,d\nu(y).$$
 Since $\tau^y (\Hla(m_j))(x)=\tau^x (\Hla(m_j))(y),$ proceeding as in the first part of the proof of \eqref{horm} we can easily check that $\sum_{j>- 2\log_2 r}|K_j(x,y)|$ is integrable over $\{y\in X\,:\,|x-y|>r\}.$
 Hence, using the dominated convergence theorem (recall that $f\in L^{\infty}$),  \begin{equation}\label{kern}\sum_{j>-2\log_2r}\,\int_{X} K_j(x,y)f(y)\,d\nu(y)=\int_{X} \sum_{j>-2\log_2r}\,K_j(x,y)f(y)\,d\nu(y).
 \end{equation}
 From \eqref{convo} it follows that
 \begin{equation}\label{multi}\mathcal T_{m_j}f(x)=H(m_j)\natural f(x)=\int_{X} K_j(x,y)f(y)\,d\nu(y),\end{equation}
 with $\mathcal T_{m_j}$ defined as in \eqref{spHan}.
 Since the Hankel transform is an $L^2(X)$-isometry, from the dominated convergence theorem we conclude  that $\sum_{j>-2\log_2r} \mathcal T_{m_j}f=\mathcal T_{m_{\infty}} f,$ where the sum converges in $L^2(X)$ and $m_{\infty}=\sum_{j>-2\log_2r}m_j.$ Hence, combining \eqref{kern} and \eqref{multi}, we obtain
 $$\mathcal T_{m_{\infty}}f(x)=\int_{X} \sum_{j>-2\log_2r}\,K_j(x,y)f(y)\,d\nu(y),$$
for $\textrm{a.e.}\,\,x$ outside
$\supp \, f.$ The function $m_{0}=m-m_{\infty}$ is bounded and compactly supported. Consequently, from \eqref{convo} we get $\mathcal T_{m_0}f(x)=\Hla(m_0) \natural f (x).$ Moreover, we see that $\sum_{j\leq -2 \log_2 r}|m_j(\la)|\leq C |m(\la)|\leq C$. Hence, from \eqref{translaform} we conclude $$\tau^{y}(m_0)(x)=\sum_{j\leq -2 \log_2 r}\tau^y(m_j)(x),$$
 so that
$$\mathcal T_{m_0}f(x)=\int_{X} \sum_{j\leq-2\log_2r}\,K_j(x,y)f(y)\,d\nu(y).$$
Then $\mathcal T_{m}f(x)=\mathcal T_{m_0}f(x)+\mathcal T_{m_{\infty}}f(x)=\int_{X} \,K(x,y)f(y)\,d\nu(y),$ as desired. \ \ $\square$
\vskip 1em

Let us finally comment that the proof of Remark \ref{remark1} goes in the same way as that of Theorem \ref{Thmoo}. The only difference is that we use Remark \ref{remark2} instead of Lemma \ref{Linterr}.

\section{Proof of Theorem \ref{Thmhar}}
We shall need the maximal-function characterization of $H^1(X).$ Define the operator $\mathcal{M} f (x)=\sup_{t>0} |T_t f(x)|,$ where $T_t f(x)= \int_{(0,\infty)^d}T_t(x,y)f(y)\,d\nu(y).$ Then we have the following proposition.
\begin{pro}
\label{Promaxat}
There exists $C>0$ such that
\begin{equation}
\label{inmaxat}
C^{-1}\|f\|_{H^1(X)}\leq \|\mathcal{M}f\|_{L^1(X)}\leq C \|f\|_{H^1(X)}.
\end{equation}
\end{pro}
\noindent The reader who is convinced that Proposition \ref{Promaxat} is true may skip Lemmata \ref{lemgaus} and \ref{lemuchi} and continue with the proof of Theorem \ref{Thmhar} on page 13. To prove the proposition we need two lemmata.
\begin{lem}
\label{lemgaus}
The heat kernel $T_t(x,y)$ satisfies the Gaussian bounds:
\begin{equation}
\label{gausbound} 0\le T_t(x,y)\leq \frac{C}{\nu (B(x,\sqrt{t}))}\exp(-c|x-y|^2\slash t),
\end{equation}
and the following Lipschitz-type estimates:
\begin{equation}
\label{heatlipsch}|T_t(x,y)-T_t(x,y')|\leq \Big(\frac{|y-y'|}{\sqrt{t}}\Big)^\delta \frac{C}{\nu(B(x,\sqrt{t}))}\exp(-c |x-y|^2\slash t), \qquad 2|y-y'|\leq |x-y|,
\end{equation}
\begin{equation}
\label{heatlipschngauss}|T_t(x,y)-T_t(x,y')|\leq \Big(\frac{|y-y'|}{\sqrt{t}}\Big)^\delta \frac{C}{\nu(B(x,\sqrt{t}))}.
\end{equation}
\end{lem}

\textbf{Proof.} Clearly, since the product of Gaussian kernels is Gaussian and $\nu$ is a product of doubling measures, it suffices to focus on $d=1.$ It is known that for $\alpha>-1/2$
\begin{equation*}
\begin{split}
T_t(x,y)& =ct^{-1}\exp(-(x^2+y^2)\slash 4t) (xy)^{-(\alpha-1)/2} I_{(\alpha -1)/2} (xy\slash 2t)\\
& =ct^{-1} \exp(-|x-y|^2\slash 4t) \exp(-xy\slash 2t)(xy)^{-(\alpha-1)/2} I_{(\alpha -1)/2}(xy\slash 2t),
\end{split}
\end{equation*}
where $I_{\mu}$ is the modified Bessel function of order $\mu.$ Using the asymptotics for $I_{\mu},$ (see \cite{3}) it is easy to see that
\begin{equation}
\label{heatest} T_t(x,y)\sim \begin{cases} t^{-(\alpha+1)\slash 2}\exp(-(x^2+y^2)\slash 4t) \ \ & \text{if } xy<t,\\
t^{-1\slash 2} (xy)^{-\alpha\slash 2} \exp(-|x-y|^2\slash 4t) \  \ & \text{if } xy\geq t.\\
\end{cases}
\end{equation}
Now, \eqref{gausbound} is a consequence of \eqref{heatest}. To prove \eqref{heatlipsch} and \eqref{heatlipschngauss}, using the identity $(x^{-\mu}I_\mu(x))'=x^{-\mu}I_{\mu+1}(x)$ and the asymptotics for $I_{\mu}$ we check that
\begin{equation*}
|\partial_y T_t(x,y)|\lesssim  \begin{cases} t^{-(\alpha+3)\slash 2}(x+y)\exp(-(x^2+y^2)\slash 4t) \ \ & \text{if } xy<t,\\
\{t^{-3\slash 2} |x-y|+t^{-1\slash 2}y^{-1}\}(xy)^{-\alpha\slash 2} \exp(-|x-y|^2\slash 4t) \  \ & \text{if } xy\geq t.\\
\end{cases}
\end{equation*}
From the above it is not hard to conclude that $$|\nabla_y T_t(x,y)|\leq \frac{C}{\sqrt{t}} \cdot \frac{1}{\nu (B(x,\sqrt{t}))}\exp(-c|x-y|^2\slash t).$$ The latter inequality easily implies \eqref{heatlipsch} and \eqref{heatlipschngauss}. $\square$\\
\indent Let $\rho (x,y)\ssb=
\inf\hspace{.5mm}\{\,\nu(B')\,|\,x,y\!\in\!B'\hspace{.5mm}\}.$ We have:

\noindent$\bullet$
\;$\rho (x,y)\hspace{-.5mm}\sim\hspace{-.5mm}\mu(B'(x,r_0)),$
\,where \,$r_0\hspace{-.75mm}=|x-y|$\ssf,
\vspace{1mm}

\noindent$\bullet$
\,$\rho(x,y)\leq A( \rho(x,z)+\rho(z,y))$
\vspace{1mm}

\noindent$\bullet$
\,$\nu(B_{\ssb \rho}^{\ssf\prime}(x,r))\sim r,$\\
\vspace{1mm}
i.e., the triple $((0,\infty)^d,d\nu,\rho)$ is a space of homogenous type.
\begin{lem}
\label{lemuchi}
Let $ K_r(x,y)=T_{t(x,r)}(x,y),$ where $t=t(x,r)$ is defined by $\nu(B(x,\sqrt t))=r.$ Then the kernel $r\,K_r$ satisfies the assumption of Uchiyama's Theorem, see \cite[Corollary 1']{Uchi}, i.e., there are constants $A,\gamma>0$ such that

\begin{equation}
\label{uchiz} K_r(x,x)\geq A^{-1}r^{-1}>0,
\end{equation}

\begin{equation}
\label{uchio} 0\le K_r(x,y)\le Cr^{-1}\bigg(\ssf1+\frac{\rho (x,y)}r\ssf\bigg)^{-1-\gamma},
\end{equation}
and
\begin{equation}
\label{uchit}\bigl|\hspace{.5mm}K_r(x,y)\ssb-\ssb K_r(x,y')\hspace{.3mm}\bigr|
\le \frac{C}{r}\bigg(\ssf1\ssb+\ssb\frac{\rho (x,y)}r\ssf\bigg)^{-1-2\gamma}
\bigg(\ssf\frac{\rho (y,y')}r\ssf\bigg)^\gamma, \qquad  \rho(y,y')\!\le\!\frac{r\ssf+\ssf \rho (x,y)}{4\ssf A}.
\end{equation}
\end{lem}

\textbf{Proof (sketch).} The inequality  \eqref{uchiz} is obvious, once we recall \eqref{heatest}.
To prove \eqref{uchio} and \eqref{uchit} we use Lemma \ref{lemgaus}. From \eqref{gausbound} we have
\begin{equation*}\begin{split}
K_r(x,y)\leq Cr^{-1}\exp(-c|x-y|^{2}\slash t).
\end{split}
\end{equation*}
Now, since
\begin{equation} \label{expa} \Big(1+\frac{\rho(x,y)}{r}\Big) \leq C\Big(1+ \frac{\nu(B(x,|x-y|))}{\nu(B(x,\sqrt{t}))}\Big)
\leq C\Big(1+\frac{|x-y|}{\sqrt{t}}\Big)^n\leq C_\varepsilon \exp(\varepsilon |x-y|^2\slash t),
\end{equation}
we get \eqref{uchio}. Observe that there is $q>0,$ such that
\begin{equation}
\label{balls} R^{q}\nu(B(x,t))\leq {C}\nu(B(x,Rt)), \qquad t>0, \quad R\geq 1.
\end{equation}
Note that we can take $q=1,$ if $\alpha_k\geq 0,$ $k=1,\ldots,d.$
The estimate \eqref{uchit} for $\rho(y,y')\geq r/(2A)$ is a simple consequence of \eqref{uchio}.
In the opposite case, i.e., $\rho(y,y')< r/(2A),$
 we first  note that \eqref{balls} implies
\begin{equation}\begin{split}\label{balls2} &\frac{\rho(y,y')}{r}  \sim \frac{\nu(B(y,|y-y'|)}{\nu(B(x,\sqrt{t}))} =
\frac{\nu(B(y,|y-y'|))}{\nu(B(y,\sqrt{t})} \cdot \frac{\nu(B(y,\sqrt{t}))}{\nu(B(x, \sqrt{t}))}
\\
& \gtrsim \Big(\frac{|y-y'|}{\sqrt{t}}\Big)^{\kappa}
\cdot \frac{\nu(B(y,\sqrt{t}))}{\nu(B(y, \sqrt{t}+|x-y|))}\gtrsim \Big(\frac{|y-y'|}{\sqrt{t}}\Big)^{\kappa}
\cdot \Big(\frac{\sqrt{t}}{\sqrt{t}+|x-y|}\Big)^{Q+d},\\
\end{split}\end{equation}
where $\kappa=q,$ if $|y-y'|\geq \sqrt{t},$ and $\kappa=Q+d,$ in the other case.
 Then \eqref{uchit} can be deduced from \eqref{heatlipsch}, \eqref{heatlipschngauss}, and \eqref{balls2}.
$\square$\\

\textbf{Proof of Proposition \ref{Promaxat}.} Since $\nu(B(x,\sqrt{t}))$ is an increasing continuous function of $t$ taking values in $(0,\infty),$ the maximal function
$$K^{\ast}f(x)=\sup_{r>0}\left|\int_{(0,\infty)^d} K_r(x,y)f(y)\,d\nu(y)\right|$$
coincides with $\mathcal{M} f.$ Now, using Lemma \ref{lemuchi} together with Uchiyama's theorem, \cite[Corollary 1']{Uchi}, we obtain a variant of the equivalence \eqref{inmaxat}, with respect to atoms corresponding to the metric $\rho.$ A simple observation that $$B\left(x,\sqrt{t(x,r)}\right)\subset B_{\rho}(x,r)\subset B\left(x,C\sqrt{t(x,r)}\right),$$ for some $C>0,$ finishes the proof. $\square$ \\
\indent Having Proposition \ref{Promaxat} we turn to prove Theorem \ref{Thmhar}.

\textbf{Proof of Theorem \ref{Thmhar}.} The proof takes some ideas from the one-dimensional case, see \cite{2}. Since the operator $\mathcal{T}_m$ maps continuously $H^1(X)$ into $\mathcal{D}'((0,\infty)^d),$ it suffices to prove that there exists a constant $C>0,$ such that for every atom $a\in H^1(X),$ we have
\begin{equation}
\label{atbound}
\|\mathcal{M}(\mathcal{T}_m a)\|_{L^1(X)}\leq C.
\end{equation}
\noindent If $a$ is an atom associated with a ball $B(y_0,r),$ then clearly,
\begin{equation}\begin{split}\label{babound}
\|\mathcal{M}(\mathcal{T}_m a)\|_{L^1(B(y_0,2r),d\nu)}&\leq \nu(B(y_0,2r))^{1/2}\|\mathcal{M}(\mathcal{T}_m a)\|_{L^2(B(y_0,2r),d\nu)}\\
  &\leq\nu(B(y_0,2r))^{1/2}\|a\|_{L^2(X)}\leq C.
\end{split}\end{equation}
\indent Fix a $C_c^{\infty}(A_{1/2,2})$ function $\psi$ satisfying
\begin{equation}\label{sumahar}
\sum_{j\in\mathbb{Z}} \psi^2(2^{-j}\la)=1, \qquad \lambda\in \mathbb R^d\backslash \{0\}.
\end{equation}
Analogously as in Section 3 we define $$m_j(\la)=\psi^2(2^{-j}(\la_1^2,...,\la_d^2))m(\la)=(\psi^2(2^{-j}\cdot)n(\cdot))(\la_1^2,...,\la_d^2).$$ In view of \eqref{babound} it is enough to show that
\begin{equation}
\label{mainhar}
\sum_{j\in\mathbb{Z}} \|\mathcal{M}(\mathcal{T}_{m_j} a)\|_{L^1((B(y_0,2r))^c,d\nu)}\leq C.
\end{equation}
Let
\begin{align*}
&m_{(j,t)}(\la)=m_j(\la)e^{-t|\la|^2}, \qquad \tilde{m}_{(j,t)}(\la)=m_{(j,t)}(2^{j/2}\la),\\
&M_{(j,t)}(x)=\Hla(m_{(j,t)})(x), \qquad \tilde{M}_{(j,t)}(x)=\Hla(\tilde{m}_{(j,t)})(x).
\end{align*}
Clearly, $M_{(j,t)}(x,y)=\tau^{y}M_{(j,t)}(x)$ are the integral kernels of the operators $\mathcal{T}_{e^{-t|\la|^2}m_j(\la)}.$ Also,
\begin{equation}
\label{eqMjt}
M_{(j,t)}(x)=(\tilde{M}_{(j,t)})_{2^{j/2}}(x),\qquad M_{(j,t)}(x,y)=2^{jQ/2}\tilde{M}_{(j,t)}(2^{j/2}x,2^{j/2}y).
\end{equation}
The following are the key estimates in the proof of \eqref{mainhar}.
\begin{lem}
\label{Lest}
There exist $\delta>0$ and $C>0$ such that for all $j\in\mathbb{Z}$ and all $r>0$ we have
\begin{equation}
\label{estM}
\int_{|x-y|>r}\, \sup_{t>0}|M_{j,t}(x,y)|\,d\nu (x)\leq C (2^{j/2}r)^{-\delta},
\end{equation}
\begin{equation}
\label{diagM}
\int_{(0,\infty)^d}\, \sup_{t>0}|M_{(j,t)}(x,y)-M_{(j,t)}(x,y')|\,d\nu (x)\leq C 2^{j/2}\,|y-y'|.
\end{equation}
\end{lem}

\textbf{Proof.} Denote
\begin{align*}
&\psi_{(j,t)}(\la)=\psi(2^{-j}(\la_1^2,\ldots,\la_d^2))e^{-t|\la|^2}, \qquad \tilde{\psi}_{(j,t)}(\la)=\psi_{(j,t)}(2^{j/2}\la),\\
&\zeta_j(\la)=\psi(2^{-j}(\la_1^2,...,\la_d^2))m(\la)=(\psi(2^{-j}\cdot )n(\cdot))(\la_1^2,\ldots,\la_d^2), \\ &\tilde{\zeta}_{j}(\la)=\zeta_j(2^{j/2}\la)=(\psi(\cdot)n(2^j\cdot))(\la_1^2,\ldots,\la_d^2).
\end{align*}
Let $\tilde{Z}_j(x)=\Hla(\tilde{\zeta}_j)(x),$ $\tilde{\Psi}_{(j,t)}(x)=\Hla(\tilde{\psi}_{(j,t)})(x).$ Arguing as in \eqref{weight}, we see that
\begin{equation}\label{weightZ}\sup_{j\in\mathbb{Z}}\|\tilde{Z}_j w^{\delta}\|_{L^1(X)}\leq C,\end{equation}
for sufficiently small $\delta>0.$
Observe that $\tilde{\psi}_{(j,t)}=n_{(j,t)},$ for some $C_c^{\infty}$ function $n_{(j,t)}$ with $\supp\, n_{(j,t)}\subset A_{1/2,2}.$ Moreover, we can check that $\sup_{(j,t)}\|n_{(j,t)}\|_{C^N}\leq C_N,$ for every $N\in\mathbb{N}.$ Hence, using \eqref{cksob} we see that for every $N>0,$ there exists $C_N^{'}$ such that
\begin{equation*}
\sup_{(j,t)}|\tilde{\Psi}_{(j,t)}(x)|\leq C_N^{'} w^{-N}(x).
\end{equation*}
From the above we see that $$|\tilde{M}_{(j,t)}(x)|=|\tilde{\Psi}_{(j,t)}\natural\tilde{Z}_j(x)|\leq C_N w^{-N}\natural|\tilde{Z}_j| (x).$$ Hence, using \eqref{weightZ} and Lemma \ref{Lemhar} we arrive at
\begin{equation*}
\int_{(0,\infty)^d}\, \sup_{t>0}|\tilde{M}_{j,t}(x,y)|w^{\delta}\,d\nu (x)\leq C.
\end{equation*}
Combining the above, together with \eqref{eqMjt} and Lemma \ref{Lfarr}, we get \eqref{estM}.\\
\indent We now turn to the proof of \eqref{diagM}. Let $\tilde{l}_{(j,t)}(\la)=e^{-t2^j |\la|^2}\psi(\la_1^2,\ldots,\la_d^2)e^{|\la|^2}$ and define $\tilde{L}_{(j,t)}(x)=\Hla(\tilde{l}_{(j,t)})(x).$ Clearly,
\begin{equation}
\label{mform}
\tilde{m}_{(j,t)}(\la)=\tilde{l}_{(j,t)}(\la)\tilde{\zeta}_{j}(\la)e^{-|\la|^2}.
\end{equation}
An argument analogous to the one presented in the previous paragraph shows that
\begin{equation*}
\sup_{j\in\mathbb{Z},t>0}|\tilde{L}_{(j,t)}(x)|\leq C_{N}^{'} w^{-N}(x).
\end{equation*}
As a consequence, there is $C>0,$ such that for every $j$
\begin{equation}
\label{lbound}
\|\sup_{t>0}|\tilde{L}_{(j,t)}\natural\tilde{Z}_j |\|_{L^1(X)}\leq C.
\end{equation}
Recalling \eqref{mform}, we obtain
\begin{equation}\begin{split}\label{heatlip}
\sup&_{t>0}|\tilde{M}_{(j,t)}(x,y)-\tilde{M}_{(j,t)}(x,y')|
\\
&=\sup_{t>0}\left|\int_{(0,\infty)^d}\,\tau^x (\tilde{L}_{(j,t)}\natural \tilde{Z}_j)(z)(T_1(z,y)-T_1(z,y'))\,d\nu(z)\right|
\\
 &\leq \int_{(0,\infty)^d}\, \tau^z \left(\sup_{t>0}|\tilde{L}_{(j,t)}\natural \tilde{Z}_j|\right)(x)|T_1(z,y)-T_1(z,y')|\,d\nu(z)|.
\end{split}\end{equation}
From \eqref{young} together with \eqref{lbound}, \eqref{heatlip} and Lemma \ref{Lneaa}, we obtain
\begin{equation}
\label{diagMt}
\int_{(0,\infty)^d}\, \sup_{t>0}|\tilde{M}_{(j,t)}(x,y)-\tilde{M}_{(j,t)}(x,y')|\,d\nu (x)\leq C |y-y'|.
\end{equation}
Now, \eqref{diagM} is a consequence of \eqref{eqMjt} and \eqref{diagMt}. $\square$\\
\indent Using Lemma \ref{Lest} and some standard arguments, as in the final stage of the proof of \cite[eq. (3.3)]{2}, we easily justify \eqref{mainhar}. Hence the proof is complete. $\square$

\end{document}